\documentstyle[12pt]{article}
\begin{document}
\begin{center}

{\bf MAXIMAL INEQUALITIES IN BILATERAL }

\vspace{3mm}

{\bf GRAND LEBESQUE SPACES OVER }\\

\vspace{3mm}

  {\bf UNBOUNDED MEASURE }\\

\vspace{3mm}

{\sc By Ostrovsky E., Rogover E.}\\

\vspace{3mm}
{\it Department of Mathematics and Statistics, Bar-Ilan University,
59200, Ramat Gan, Israel.}\\
e-mail: \ galo@list.ru \\

{\it Department of Mathematics and Statistics, Bar-Ilan University,
59200, Ramat Gan, Israel.}\\
e - mail: \ rogovee@gmail.com \\

\vspace{4mm}

{\it Abstract}\\

\end{center}

\vspace{3mm}

  In this paper non-asymptotic exact rearrangement invariant norm
estimates are derived for the maximum distribution  of  the family
elements of some rearrangement invariant (r.i.) space over unbounded
measure in the entropy terms and in the terms of generic chaining.\par
 We consider some applications in the martingale theory and in the
theory of Fourier series.\par

\vspace{3mm}


\vspace{3mm}


\vspace{4mm}

{\bf 1. Introduction. Notations. Statement of problem.} \par

 Let $ (X, \Sigma,\mu) $ be a measurable space with non-trivial
measure $ \mu: \ \exists A \in \Sigma, \mu(A) \in (0,\mu(X)). $ \par
 We will assume that $ \mu(X) = \infty $
and that the measure $ \mu \  $ is $ \sigma - $ finite and diffuse:
$ \forall A \in \Sigma, 0 < \mu(A) < \infty \ \exists B \subset A, \mu(B) =
\mu(A)/2. $ \par

 Let also $ T = \{t \} $ be {\it arbitrary } set and $ Y = Y(t,x) = Y(t) $
be some function of a variables $ t $ and $ x $  such that for all the values
$ t \in T $ the function $ Y =  Y(t, x) $ is measurable as a function on $ x $
and is separable. \par
{\bf Definition 1.} The function $ Y = Y(t,x) $ is called separable relatively
the  variable $ t; \ t \in T, $ if there exists a countable
subset $ \tilde{T} $ of a set $ T: \  \tilde{T} = \{ t_1, t_2, \ldots \}
\subset T $ such that for arbitrary closed set $ Q $ on the space  $ R = R^1 $

$$
\cap_{t \in \tilde{T}} \{ x: Y(t) \in Q \} \sim \cap_{t \in T}
\{ x: Y(t) \in Q \}. \eqno(1.1)
$$

 Here and further the set equivalence $ A \sim B, \ A,B \subset X $ means that
both  the sets $ A $ and $ B $ are measurable: $ A \in \Sigma, \ B \in \Sigma $
and

$$
 \mu \{ (A \setminus B) \cup ( B \setminus A) \} = 0.
$$

 As a rule, the set  $ \tilde{T} $ is enumerable dense subset of $ T $ relatively  some  distance (or semi-distance)  $ r = r(t,s), \ t,s \in T $
on the set $ T. $ We will  call in this case the subset $ \tilde{T} $ the
{\it separante } of the set $ T $ and will write

$$
 \tilde{T} = sep(T, r).  \eqno(1.2)
$$

 For example, if the function $ Y = Y(t, \cdot) $ is continuous a.e. on the
variable $ t $ in the distance $ r, $ and the metric space $ (T,r) $ is separable, then $ Y(t,x) $ is separable. Further, if the set $ T $ is
the union of some sequence subsets $ S_m, m = 1,2,\ldots, M, \ M \le \infty $
of the set $ T $

$$
T = \cup_{m=1}^M S_m
$$

 and the function $ Y $ is separable on the sets $ S_m, $ then $ Y = Y(t,x) $
is separable on the set $ T. $ \par

 Let us define

$$
\overline{Y} = \overline{Y}(x) = \sup_{t \in T} Y(t,x).  \eqno(1.3)
$$

 It is easy to verify, as  in the theory of random processes,
 that if the function $ Y = Y(t,x) $ is separable, then
$ \overline{Y}(x) $ is measurable function on the variable $ x. $ \par

\vspace{2mm}

{\sc Further we will assume that our function } $ Y = Y(t,x) $ {\sc is
separable over some suitable dense set } $ \tilde{T}. $ \par

\vspace{2mm}

 Let also $ G $ be some rearrangement invariant (r.i.) space with a
norm $ ||\cdot||G $ over our triplet $ (X, \Sigma,\mu), $ for instance,
$ L_p = L_p(X, \Sigma,\mu), $ Orlicz, Marcinkiewicz, Lorentz or
Grang Lebesque spaces etc. \par

{\bf Our aim is  obtaining the $  G \ - $ norm estimation for
$ \overline{Y}: \ || \ \overline{Y} \ ||G $ through some simple
rearrangement invariant parameters of source function } $ Y(t,x). $ \par

  The important examples of these parameters are: the maximal value

$$
\sigma = \sigma(G) \stackrel{def}{=} \sup_{t \in T} ||Y(t, \cdot)||G
\eqno(1.4)
$$
and  the so-called $ G \ - $ distance (more exactly, semi-distance) $ d_G = d_G(t,s) $ on the set $ T: $

$$
d_G = d_G(t,s) \stackrel{def}{=} ||\ Y(t,\cdot) \  - \ Y(s, \cdot) \ ||G.
\eqno(1.5)
$$

 Recall that the semi-distance $ d = d(t,s), \ s,t \in T $ is, by definition,
non-negative symmetrical numerical function, $ d(t,t) = 0, \ t \in T, $
satisfying the triangle inequality, but the equality $ d(t,s) = 0 $
does not means (in general case) that $ s = t. $ \par

  It is evident that if $ \sigma(G) < \infty, $ then $ d_G(t,s) \le
2 \sigma(G). $ \par

 Notice that the case $ \mu(X) = 1 $ (the "probabilistic case" ) is well
investigated in the theory of random fields, see, for example,[1], [2], [3],
[4], [5], [6], [41] etc. The obtained there results
may be used here as illustration of precision of estimations of this article.\par

We will use widely further the notion of {\it fundamental function}
$ \phi(G, \delta), \ \delta \in (0,\infty) $ of the r.i. space $ G. $ Recall
 that  by definition

$$
 \phi(G, \delta) = || \ I(A) \ ||G, \ \mu(A) = \delta
$$
and $ I(A) = I(A,x) = 1, x \in A, \ I(A) = I(A,x) = 0, x \notin A. $ \par
 This notion play a very important role in the theory of interpolation of
operators, theory of Fourier series, theory of approximation etc. See, for example, [30],  [23], [42] etc.\par
 Let the set $ T $ relatively some semi-distance $ r = r(t,s) $ be precompact
set. We denote by $ N =  N(T, r, \epsilon) $ the minimal number of closed
$ r \ - $ balls $ B(r; t_j, \epsilon), \ t_j \in T  $ with the center $ t_j $ and the radius $ \epsilon, \ \epsilon > 0: $

$$
B(r; t_j, \epsilon) = \{t, \ t \in T, \ r(t, t_j) \le \epsilon \}
$$

covering the set $ T: $

$$
N(T, r, \epsilon) \stackrel{def}{=} \min \{K, \ \exists \{t_j \}, j = 1,2,\ldots, K; \ t_j \in T, \ T \subset \cup_{j=1}^K B(r; t_j, \epsilon) \}.
$$

 The (natural) logarithm of $ N(T, r, \epsilon):  \
H(T, r, \epsilon) =  \log N(T, r, \epsilon) $ is called entropy of $ T $
in the distance $ r, $  and the value (finite or infinite)

$$
\kappa = \kappa_r(T) \stackrel{def}{=}
\overline{\lim}_{\epsilon \to 0+} \frac{ H(T, r, \epsilon)}{ |\log \epsilon|}
$$

is called the {\it dimension} $ T $ in the distance $ r:$
$$
\kappa_r(T) = dim_r(T).
$$

\vspace{3mm}

{\bf 2. Grand Lebesque spaces.} \par

\vspace{3mm}

  We define as usually for arbitrary measurable function $ f: X \to R^1 \ $

$$
{\bf E }f = \int_X f(x) \mu(dx); \ p \ge 1 \ \Rightarrow
$$

$$
|f|_p = {\bf E}^{1/p} \left( |f|^p \right)   =
\left(\int_X |f(x)|^p \ \mu(dx) \right)^{1/p};
$$
$ L_p = L(p) = L(p; X,\mu) = \{f, |f|_p < \infty \}. $ \par

 Let $ a = const \ge 1, b = const \in (a,\infty], $ and let $ \psi = \psi(p)
= \psi(p; a,b) $
be some strong positive: $ \psi(p) \ge 1 $ bounded in each open
subinterval $ (c,d), \ a < c < d < b $ logarithmical convex
on the {\it open} interval $ (a,b) $ function. \par
 We will denote the set of all such a function  by $ \Psi: \ \Psi =
\Psi(a,b) = \{ \psi \} = \{\psi(\cdot; a,b) \}. $ \par

\vspace{3mm}

{\bf Definition 2.} The space
$ BGL(\psi) = G(\psi) = G(X,\psi) = G(X,\psi, \mu) = G(X,\psi,\mu, a,b) $
(Bilateral Grand Lebesque space)
consist on all the measurable functions $ f: X \to R $ with finite norm
$$
||f||G(\psi) \stackrel{def}{=} \sup_{p \in (a,b)} \left[ |f|_p/\psi(p) \right].
\eqno(2.1)
$$
 We can define formally in the case $ a = b \in [1, \infty) \ G(\psi) = L_a.$ \par

 Suppose that there exist a {\it pair} of numbers $ (a,b), \ 1 \le a < b \le \infty $ such that

$$
\forall p \in (a,b) \ \Rightarrow \ |Y(t, \cdot)|_p < \infty
$$
 and such that

$$
\forall \Delta > 0 \ \Rightarrow \ \sup_{t \in T} |Y(t, \cdot)
|_{a \ - \ \Delta} = \infty \eqno(2.2)
$$
and
$$
\forall \Delta > 0 = \sup_{t \in T} |Y(t, \cdot) |_{b + \Delta} = \infty
 \eqno(2.3)
$$

where in the case $ a = 1 $ the condition (2.2) is absent and in the case
$ b = \infty $ the condition (2.3) is absent. \par

 Then we can define the following {\it natural choice} of a function
$ \psi_0(p) $ as follows:

$$
\psi_0(p) \stackrel{def}{=} \sup_{t \in T} |Y(t, \cdot)|_p. \eqno(2.4)
$$

 The spaces $ G(\psi), \ \psi \in U\Psi $ are non-trivial: arbitrary bounded
$ \sup_x |f(x)| < \infty $
measurable function $ f: X \to R $ with finite support: $ \mu(supp \ |f|)
< \infty $ belongs to arbitrary space $ G(\psi). $ \par
 We denote as usually $ supp \ \psi = \{p: \ |\psi(p)| < \infty \}. $ \par

  The detail investigation of these spaces see, for example, in [14], [15],
[17], [18], [34], [42]
 etc.\par
  It is knows (see [42]) that the BGL spaces in general case does not
coincide  with classical r.i. spaces: Lorentz, Marcinkiewicz, Orlicz
spaces. It is obvious that BGL spaces does not coincide with recently
appeared Grand Orlicz, modular and variable Lebesgue spaces, as long as
both the last spaces are not, in general case, rearrangement invariant
(see [37], [38], [39]). \par

  The  BGL spaces are used, for example, in the theory of probability [2]- [7],
 [8] - [10], [42]; theory of PDE [14], [15], functional analysis [11], [12],
 [42], theory of Fourier series [23], [30], theory of martingales
 [14], [15], [16] etc.\par

  If we choose as the r.i. space $ G $ the space $ G(\psi_0), $ then
$ \sigma(G(\psi_0) ) = 1; $ and we can introduce the so-called {\it natural}
distance

$$
d_0(t,s) \stackrel{def}{=} ||Y(t, \cdot) \ - \ Y(s, \cdot)||G(\psi_0), \
t,s \in T.
$$

 This approach in the probabilistic case was introduced by [36] for
Gaussian random fields; more general case was considered in [4].\par

 The fundamental function of BGL spaces may be calculated by the formula:

$$
\phi(G(\psi), \delta) = \sup_{p \in (a,b)} [ \delta^{1/p} /\psi(p)].
$$
 Many examples of $ G(\psi) $ spaces and its fundamental functions see in [42]. As a particular case $ G(\psi) $ space may coincide with arbitrary
exponential Orlitz's space.\par

 The spaces $ G(\psi, a,b) $ are non-separable  and non-reflexive
([42]), but they satisfy the Fatou property.  Namely, the following property about these spaces is true.\par

\vspace{2mm}

{\bf Proposition 1.} \ The $ G(\psi) $ space satisfies the Fatou property. \par

{\vspace{2mm}

{\bf Proof.} Recall at first that the Fatou property of some r.i. space
$ G $ over source triplet $ (X, \Sigma, \mu) $ denotes that for arbitrary
non-increasing sequence of non-negative functions $ \{  f_n \}
= \{ f_n(x), \ x \in X  \} $
belonging
to the space $ G $ and such that  as $ n \uparrow \infty $

$$
f_n(x) \uparrow f(x), \ \sup_n ||f_n||G  < \infty  \eqno(2.5)
$$
it follows

$$
||f_n||G \uparrow ||f||G. \eqno(2.6)
$$
 Let $ G = G(\psi) $ and let the sequence of measurable
functions $ \{f_n \} = \{f_n: X \to R \} $ satisfies the condition (2.5).
As long as the space $ L_p(X, \mu) $ satisfies the Fatou property, we have:

$$
\sup_n ||f_n||G(\psi) = \sup_n \sup_{p \in (a,b)} [ |f_n|_p/\psi(p) ] =
$$

$$
\sup_{p \in (a,b)} \sup_n [|f_n|_p/\psi(p) ] = \sup_{p \in (a,b)}[|f|_p
/\psi(p)] = ||f||G(\psi),
$$
Q.E.D. \\
 As a simple consequence: it follows from theorem of Mityagin  \ - \ Kalderon that
the space $ G(\psi) $ is interpolation space between spaces $ L_1(X,\mu) $
and $ L_{\infty}(X,\mu). $ See in detail [11], [12].\par

\vspace{2mm}

{\bf 3. Main results.} \par

\vspace{2mm}

{\bf A.} Generic chaining theory in our case.\par

Now we recall, modify and rewrite some definition from the "generic
chaining" theory, belonging to X.Fernique [1] and M.Talagrand [6] - [10]. \par

 Let $ (G, ||\cdot||G $ be some r.i.space over $ (X, \Sigma, \mu) $ and let

$$
\tilde{T} = sep(T, d_G).
$$

{\bf Definition 3}. The generic chaining $ W $ is, by definition, the partition
of the   set  $ \tilde{T} $ into a sequence of finite subsets $ \{ Q_k \}: $

$$
\tilde{T} = \cup_{k=0}^{\infty} Q_k,
$$

where $ |Q_k| \stackrel{def}{=} card \ (Q_k) < \infty. $ Notation:
$ W = \{ Q_k \}. $ \par
 Without loss of generality we can and will assume that $ Q_0 = \{t_0 \}, $
where

$$
\sigma(G) = \sup_{t \in T} ||Y(t, \cdot) || G
= ||Y(t_0, \cdot)||G.
$$

 For any element $ t \in T $ we denote arbitrary, but
fixed (non-random) element $ \pi_k(t)$ of a subset  $ Q_k $
 such that
$$
d_G (t, \pi_k(t) ) = \min_{s \in Q_k }d_G(t,s). \eqno(3.1)
$$

Thus,

$$
|| Y(t, \cdot)  \ - \ Y(\pi_k(t), \cdot)||G \le d_G(t, \pi_k(t)). \eqno(3.2)
$$

Let us denote for some partition $ W = \{ Q_k \} = \{ Q(k) \} $

$$
\Lambda(T, G, W) = \sum_{k=0}^{\infty} || \max_{t \in Q_k}
( Y(\pi_k(t), \cdot)- Y(\pi_{k-1}(t), \cdot) ||G.
$$

\vspace{2mm}

{\bf Proposition  2.} \par

\vspace{2mm}

$$
|| \overline{Y} ||G \le \inf_W  \Lambda(T, G,W). \eqno(3.3)
$$

{\bf Proof } is very simple.  Let $ R $ be arbitrary partition.
Since the function $ Y = Y(t,x) $ is presumed to be separable,
we have a.e.:

$$
\overline{Y} = \lim_{ M \to \infty} \max_{t \in \cup_{k=1}^M Q(k)} Y(t, x)
\le
$$

$$
\lim_{M \to \infty}\sum_{k=0}^M \max_{t \in Q_k} ( Y(\pi_k(t), x) \ - \
 Y(\pi_{k-1} (t), x)).
$$

We find
using the triangle inequality for the $ G \ - $ norm

$$
|| \overline{Y} ||G \le  \Lambda(T, G, W). \eqno(3.4)
$$

Since the partition $ W $ is arbitrary, we get to the (3.3) after the minimization
over $ W. $ \par
 Following, we need to estimate the $ G \ - $ norm for the
maximal value of finite set of a functions. At first we use the so-called
Pizier technique. \par

{\bf B.} (Finite case). We suppose here that the set $ T $ is finite: $ T =
\{t_1, t_2, \ldots, t_m \}; $  on the other words, $ card(T) = m < \infty, $
and  assume that for some $ p \in [1,\infty) $

$$
 \max_{j = 1,2,\ldots,m} |Y(t_j, \cdot)|_p < \infty.
$$

\vspace{2mm}

{\bf Proposition 3.}

\vspace{2mm}

{\it  We provide the following generalization of famous Pizier's [10] inequality:}

$$
| \overline{Y} |_p \le \max_{j = 1,2,\ldots,m} |Y(t_j, \cdot)|_p \cdot m^{1/p}.
\eqno(3.5)
$$

{\bf Proof.} Indeed, assume for simplicity $ | \ Y(t_j) \  |_p \le 1. $ We get:

$$
[\overline{Y}]^p= \max_{j=1,2, …,m} [Y(t_j, \cdot)]^p
 \le \sum_{j=1}^m [Y(t_j, \cdot)]^p;
$$

 $$
 |\overline{Y} |^p_p \le \sum_{i=1}^m | Y(t_i, \cdot)|^p_p  \le m.
 $$

\vspace{2mm}

{\bf C.} ( Generalization of finite case). \par

Let $ \psi, \zeta, \nu $ be three function from the set $ \Psi(a,b) $ such
that

$$
\zeta(p) = \psi(p) \ \nu(p), \ p \in (a,b).
$$

 We suppose again here that the set $ T $ is finite: $ T =
\{t_1, t_2, \ldots, t_m \}; $
and  assume that for some $ p \in [1,\infty) $

$$
 \max_{j = 1,2,\ldots,m} |Y(t_j, \cdot)|_p < \infty.
$$

\vspace{2mm}

{\bf Proposition 4.}

\vspace{2mm}

$$
|| \overline{Y} ||G(\zeta) \le \max_{i=1,2,…,m} ||f_i||G(\psi) \cdot
 \phi(G(\nu), m). \eqno(3.6)
$$
{\bf Proof.}  We may use the inequality (3.5), estimating the values
$ |f_i|_p $ as

$$
|f_i|_p \le ||f_i||G(\psi) \cdot \psi(p),
$$
on the basis of definition the norm in the $ G(\psi) $ space. We have:

$$
| \overline{Y} |_p \le  \max_{i=1,2,…,m} ||f_i||G(\psi) \cdot \psi(p)
\cdot m^{1/p}.
$$

 Dividing by $ \zeta(p) $ and tacking supremum  over $ p \in (a,b), $
we receive:

$$
|| \overline{Y} ||G(\zeta) \le \max_{i=1,2,…,m} ||f_i||G(\psi) \cdot
\sup_{p \in (a,b)} \frac{m^{1/p} }{\nu(p)} =
$$

$$
 \max_{i=1,2,…,m} ||f_i||G(\psi) \cdot  \phi(G(\nu), m),
$$
 Q.E.D. \par

\vspace{3mm}

 Let now and further $ \theta $ be some {\it fixed } number inside the interval
(0, 1), for example, $ \theta = 1/2 $ or $ \theta = 1/e.$ We suppose
for some $ p \ge 1 $

$$
  \sup_{t \in T} | Y(t,\cdot) |_p < \infty,
$$

and denote

$$
d_p(t,s) \stackrel{def}{=} | \ Y(t, \cdot) \  - \ Y(s, \cdot) \ |_p.
$$

 We consider here as the set $ Q_k $ and consequently the partition $ W $
in (3.3)
the minimal $ \theta^k $ set of the space $ T $ under the
distance $ d_p; $ recall that the quantity of its element is equal to
$ N(T,d_p, \theta^k). $ \par

\vspace{2mm}

{\bf Proposition 5.}

\vspace{2mm}

$$
| \ \overline{Y} \ |_p \le \sum_{k=1}^{\infty}
\theta^{k-1} N^{1/p}(T,d_p, \theta^k). \eqno(3.7)
$$

{\bf Proof.} This proposition follows immediately from proposition 2 and our
generalization of Pizier inequality (3.5):

$$
|\max_{t \in Q(k)} (Y(\pi_k(t), \cdot) \ - \ Y(\pi_{k \ - \ 1}(t), \cdot))|_p
\le \theta^{k-1} N^{1/p}(T,d_p, \theta^k)
$$
after summing  over $ k. $ \par

{\bf Remark 1.} We can rewrite the inequality (3.7) as follows:

$$
| \ \overline{Y} \ |_p \le  \inf_{\theta \in (0,1)} \sum_{k=1}^{\infty}
\theta^{k-1} N^{1/p}(T,d_p, \theta^k).
$$

\vspace{2mm}

{\bf Formulation of main result.} \par

\vspace{2mm}

 Let  as in the proposition 4
$ \ \psi, \zeta, \nu $ be three function from the set $ \Psi(a,b) $
 Fixing some pair $ a,b: \ 1 \le a < b \le \infty $ and a three functions
$ \zeta(\cdot), \psi(\cdot), \nu(\cdot) $ from the space $ \Psi(a,b) $
such that

$$
\zeta(p) = \psi(p) \ \nu(p), \ p \in (a,b),
$$
we assume that

$$
\sup_{t \in T} || \ Y(t,\cdot) \ ||G(\psi) < \infty,
$$
and denote

$$
d_{\psi}(t,s) = || \ Y(t, \cdot) - Y(s, \cdot) \ ||G(\psi).
$$

 For example, $ \psi(p) $ may coincide with  the natural function
$ \psi_0(p). $\par

 We consider in this section as the set $ Q_k $ and consequently the
partition $ W $ in (3.3)
the minimal $ \ \theta^k \ - $ set of the space $ T $ under the
distance $ d_{\psi}; $ recall that the quantity of its element is equal to
$ N(T,d_{\psi}, \theta^k). $ \par

\vspace{2mm}

\bf Theorem 1.} \par

\vspace{2mm}

$$
|| \ \overline{Y} \ ||G(\zeta) \le  \inf_{\theta \in (0,1)} \sum_{k=1}^{\infty}
\theta^{k - 1} \phi \left(G(\nu),N(T,d_{\psi}, \theta^k ) \right). \eqno(3.8)
$$

{\bf Proof } is  at the same as in the proposition 5; instead the Pizier
inequality (3.5) we use its generalization (3.6). \par

 Note that it follows from conclusion of Theorem 1 the {\it continuity}
of $ Y(t) $ with probability one in the semi-distance $ d_{\psi}: $

$$
\mu \{x:  Y(\cdot,x) \notin C \left(T,d_{\psi} \right) \} = 0;
$$
$ C(T,d) $ denotes as usually the space of all continuous with
respect to the semi-distance $ d $ functions $ f: T \to R. $ \par

 The conditions of theorem 1 in the "probabilistic" case $ \mu(X) = 1 $
are equivalent to the so-called condition of the "convergence of the
majoring integral", see [7], [8]. \par

\vspace{3mm}

{\bf Examples.} \\

\vspace{2mm}

{\bf Example 1.} Let under the conditions of theorem 1 for all values
$ \epsilon \in (0, \theta) $ and for some $ \kappa = const > 0 $

$$
N(T, d_{\psi}, \epsilon) \le C \ \epsilon^{-\kappa }. \eqno(3.9)
$$

Denote for the values $ p > \max(\kappa,1) $

$$
  \psi^{( \kappa)}(p) = \psi(p) \cdot \frac{p}{p - \kappa}. \eqno(3.10)
$$
 As long as

$$
N(T, d_p, \theta^k ) \le N(T, d_{\psi}, \theta^k/\psi(p) ),
$$

we obtain after some calculations using the result (3.7) of the proposition
5:

$$
| \overline{Y} |_p \le \psi(p) + [\psi(p)]^{ \kappa/p} \sum_{k=1}^{\infty}
\theta^{k(1 - \kappa/p)} \le
$$

$$
\psi(p) + C[\psi(p)]^{ \kappa /p} \cdot
\left[ \theta^{ \kappa /p} \ - \ \theta)^{-1} \right] \le
$$

$$
C \psi(p) \left[ 1 + ( \theta^{ \kappa /p} \ - \ \theta )^{-1} \right] \le
C \psi^{( \kappa)}(p), \ C = const. \eqno(3.11)
$$

 Therefore, under considered conditions

$$
||\overline{Y}||G \left(\psi^{( \kappa)} \right) \le C \sup_{t \in T}
|| Y(t, \cdot)||G(\psi). \eqno(3.12)
$$

Since

$$
|| \ \overline{Y} \ ||G(\psi) \ge \sup_{t \in T} || Y(t, \cdot)||G(\psi),
$$

we conclude that the estimation (3.12) is exact up to multiplicative constant
in the case if $ \psi(\cdot) \in \Psi(a,b), \ \zeta(p) = \psi(p), $
where $ \kappa < a; $ the last condition is satisfied automatically if
$ \kappa < 1. $ \par

 In the case if for all values $ \epsilon < \theta $

$$
N(T, d_{\psi}, \epsilon) \le C \ \epsilon^{-\kappa(1)} \
|\log \epsilon |^{- \kappa(2) }, \eqno(3.13)
$$

$ \kappa(1) = const > 0, \kappa(2) = const < \kappa(1), $ we obtain after
some calculations denoting for the values $ p > \kappa(1), \ p \in (a,b) $

$$
\psi_{\kappa(1), \kappa(2)}(p)=  \left[ \frac{p}{p \ - \ \kappa(1)} \right]^
{1 - \kappa(2)/\kappa(1)} \cdot \psi(p):
$$

$$
|| \overline{Y} ||G \left(\psi_{\kappa(1), \kappa(2)} \right) \le  C \ \sup_{t \in T}
|| Y(t, \cdot)||G(\psi). \eqno(3.14)
$$
 In the case if the condition (3.13) is satisfied and
$ \kappa(1) = const > 0, \kappa(2) = \kappa(1), $ we conclude denoting

$$
\psi_{l, \kappa(1), \kappa(2)}(p)=  \left| \frac{ \log (p \ - \ \kappa(1))}
{ \log(p) } \right|_+ \cdot \psi(p),
$$

$ z_+ = \max(z,1): $

$$
|| \overline{Y}||G \left(\psi_{l, \kappa(1), \kappa(2)} \right) \le  C \ \sup_{t \in T}
|| Y(t, \cdot)||G(\psi). \eqno(3.15)
$$

 Finally, in the case if the condition (3.13) is satisfied and
$ \kappa(1) = const > 0, \kappa(2) > \kappa(1), $ we conclude:

$$
||\overline{Y}||G(\psi) \le C \  \sup_{t \in T}
|| Y(t, \cdot)||G(\psi). \eqno(3.16)
$$

 The estimations (3.14), (3.15), (3.16) it follow from Theorem 1 and the
elementary inequalities (3.17.1), (3.17.2), (3.17.3),
where we denote

$$
S_{\beta}(q) = \sum_{k=1}^{\infty} q^k \ k^{\beta}, \ q \in [1/2, 1), \
\beta = const:
$$

$$
\beta > - 1 \ \Rightarrow \ S_{\beta}(q) \le C(\beta) \ (1 \ - \ q)^{- 1 - \beta};
\eqno(3.17.1)
$$

$$
\beta = - 1 \ \Rightarrow \ S_{\beta}(q) \le C \ |\log (1 \ - \ q)|;
\eqno(3.17.2)
$$

$$
\beta < - 1 \ \Rightarrow \ S_{\beta}(q) \le C(\beta). \eqno(3.17.3)
$$

\vspace{2mm}

{\bf Example 2.} Exponential Orlicz spaces. \par
 We consider here as a space $ G $ a so-called exponential Orlicz spaces.\par

\vspace{2mm}

{\bf Definition 3.} We introduce the $ N(a, \beta) = N(a, \beta; u),
a \ge 1, \beta > 0 $ as an Orlicz's function such that

$$
u \to 0 \ \Rightarrow \ N(a, \beta; u) \sim C_1 |u|^a;
$$

$$
|u| \to  \infty \ \Rightarrow N(a, \beta; u) =
\exp \left(C_2 |u|^{1/\beta} \right).
$$

 The correspondent Orlicz space defined over source triple with $ N \ - $
Orlicz function $ \Phi(u) = \Phi(a,\beta; u) $ will denoted  as $ Or(a, \beta) $
and the norm of a (measurable) function $ f: X \to R $ in this space
will denoted as

$$
||f||G(a, \beta) = ||f||Or(a, \beta) = ||f||Or(\Phi(a,\beta; \cdot). \eqno(3.18)
$$

 Let $ a = const \ge 1,  \beta(1), \beta(2)= const, \ 0 < \beta(1) < \beta(2)
< \infty. $ Suppose that

$$
\sup_{t \in T} || Y(t, \cdot) ||G(a, \beta(1) ) < \infty
$$
and introduce a distance $ d_{a, \beta(1)} (t,s) $ by the formula

$$
d_{a, \beta(1)} =
d_{a, \beta(1)} (t,s) = || Y(t, \cdot) \ - \ Y(s, \cdot)|| G(a, \beta(1) ).
$$

We assert:

$$
 || \overline{Y} || G(a, \beta(2)) \le
 C \ \sup_{t \in T} || Y(t, \cdot) ||G(a, \beta(1) ) \ \times
$$
$$
\inf_{\theta \in (0,1)} \sum_{k=1}^{\infty} \theta^{k \ - \ 1}
H^{\beta(2) - \beta(1)} \left(T, d_{a, \beta(1)},\theta^k \right). \eqno(3.19)
$$
 Recall that $ H(T,d,\epsilon) = \log N(T,d, \epsilon). $ \par
The proof of (3.19) it follows from theorem 1 and  from the fact that the
space $ Qr(a, \beta) $ coincides up to the norm equivalence with some

$ G(\psi) = G(\psi; a, \infty) $ space:

$$
||f||G(a, \beta) = ||f||Or(a, \beta) \asymp \sup_{p \ge a}
\frac{|f|_p}{p^{\beta} }.
$$
 See for example [23], [42] where is formulated and proved more general assertion. \par

 Note that the inequality (3.19) is alike to the famous Dudley condition for
continuity of  Gaussian random field [36]. \par

Note also that the condition

$$
\inf_{\theta \in (0,1)} \sum_{k=1}^{\infty} \theta^{k \ - 1 \ }
H^{\beta(2) - \beta(1)} \left(T, d_{a, \beta(1)},\theta^k \right) < \infty
\eqno(3.20)
$$

is  satisfied if for example
$$
dim (d_{a, \beta(1)}, T) < \infty.
$$

\vspace{3mm}

{\bf 4. Generalization on the moment rearrangement spaces.}\par

\vspace{3mm}

  Let $ (G, \ ||\cdot||G) $ be some r.i. space defined over our triplet
$ (X, \Sigma, \mu). $  We reproduce in this section the notion of the so-called
{\it moment rearrangement invariant (m.r.i.) } space  from [29] and
consider the generalization of maximal inequality on m.r.i. spaces.\par

\vspace{2mm}

 {\bf Definition 4.} \par

\vspace{2mm}	

  We will say that the r.i. space $ G = G(m) = G_m $ with the norm
$ ||\cdot||G = ||\cdot||G(m) $ is
 moment rearrangement invariant space, briefly: m.r.i. space, or \\
$ G = G(m) = (G, \ ||\cdot||G) \in m.r.i., $
 if there exist a real constants $ a, b; 1 \le a < b \le \infty, $ and some
 rearrangement invariant norm  $ < \cdot > $ defined on the space of a real functions defined on the interval $ (a,b), $ not necessary to be finite on all the functions, such that

  $$
  \forall f \in G \ \Rightarrow || f ||G = < \ h(\cdot) \ >, \ h(p) = |f|_p. \eqno(4.1)
  $$

    We will write for considered m.r.i.  spaces $ (G, \ ||\cdot||G) $

   $$
       (a, b) \stackrel{def}{=} supp(G),
   $$
   "moment support"; not necessary to be uniquely defined.

   There are many r.i. spaces satisfied the condition (4.1) aside from
$ G(\psi) $ spaces:
 exponential Orlicz's spaces, Marcinkiewicz spaces, interpolation spaces (see [29], [33], [35]).\par

    In the article [32] are introduced the so-called $ Q(p,\alpha) $ spaces
  consisted on all the measurable function $ f: T \to R $  with finite norm
   $$
   ||f||_{p,\alpha} = \left[ \int_1^{\infty} \left(\frac{|f|_x}{x^{\alpha}}
    \right)^p \ \nu(dx) \right]^{1/p},
   $$

   where $ \nu $ is some Borelian measure.\par
    Astashkin in [33] proved that the space $ Q(p,\alpha) $ in the case
   $ T = [0,1] $ and $ \nu = m, \ m $ is Lebesgue measure coincides with the
   Lorentz $ \Lambda_p( \log^{1-p \alpha}(2/s) ) $ space.  Therefore,
     both this spaces are m.r.i. spaces.\par
  Since for arbitrary real-valued continuous function $ f $ defined on the set
 $ [0,1] $

  $$
 ||f||C[0,1] = \sup_{t \in [0,1]}|f(t)| = \lim_{p \to \infty } |f|_p = \sup_{p \in [1,\infty) } \ | \ |f|_p \ |,
  $$
  the space $ C[0,1] $ is m.r.i. space with $ supp( C[0,1]) = [1,\infty) $
 or equally, e.g., $ supp (C[0,1]) = [3, \infty). $ \par

 But there exist rearrangement invariant spaces without m.r.i. property
[29]. \par

  Let $ G = G_m $ be some m.r.i. space and suppose for all values $ p \in (a,b) $
$$
\sup_{t \in T} | Y(t, \cdot)|_p < \infty.
$$
Denote as in the section 3

$$
d_p(t,s) = | \ Y(t,\cdot) \ - \ Y(s, \cdot) |_p.
$$

{\bf Proposition 6.}

 We denote also
$$
g(p) =  \inf_{\theta \in (0,1)} \sum_{k=1}^{\infty}
\theta^{k-1} N^{1/p}(T,d_p, \theta^k).
$$

 It follows from the definition of m.r.i. spaces (4.1) and from the
proposition 5 that

$$
|| \ \overline{Y} \ ||G(m) \ \le \ <g>. \eqno(4.2)
$$

\vspace{2mm}

{\bf 5.  Application to the martingale theory over the spaces with infinite measure.} \par

\vspace{2mm}

 Let $ (S_n, F_n) = (S(n), F(n)) $ be a martingale, i.e. a monotonically
 non – decreasing sequence of $ F_n \ - $ sigma - subalgebras of $ \Sigma $
and  $ F_n  = F(n) $ measurable functions $ S_n $ such that $ {\bf E} S_{n+1}/F_n = S_n \ $ {\bf a.e.}.  \par
 We define formally $ S(0) = S_0 = 0; \ F(0) = F_0 = \{\emptyset, X \}. $ \par

 In this section we will use also the probabilistic notations

$$
{\bf Var} \ f = {\bf Var} ( f ) =  {\bf E}( f \  - \ {\bf E} f )^2 =
|f  \ - \ {\bf E} f|_2^2
$$

and notation $ {\bf E } f /F $  for the conditional expectation. \par
 Denote

$$
\sigma(n) = \left[ {\bf Var} (S_n)  \right]^{1/2}
$$
 and suppose the function $ n \to \sigma(n) $ be {\it regular} varying:

$$
\sigma(n) = n^{\gamma} L(n), \gamma = const > 0,
$$
where $ L = L(n) $ is {\it  slowly} varying as $ n \to \infty: $
$$
\forall C > 0 \ \Rightarrow  \ \lim_{n \to \infty} L(Cn)/L(n) = 1.
$$
 It is obvious that

$$
\sigma^2(n) = \sum_{k=1}^n || S(k) \ - \ S( k \ - \ 1) ||^2_2.
$$

 The $ L_p \ - $  theory of conditional expectations  and theory of
martingales in the case
$ \mu(X) = \infty $ and some its applications see, for example, in the
book  [24], pp. 330 - 347; see also [25], [26]. \par
 The Orlicz's norm estimates for martingales are used in the modern
non - parametrical statistics, for example, in the so - called regression
problem ([4], [42] etc). \par
 We recall here the famous inequality of Doob:

$$
p > 1 \ \Rightarrow \left| \sup_{ n \in [1, N] } |S_n| \  \right|_p \le
\sup_{n \in [1, N]} \left[  |S_n|_p  \ p/(p-1) \right], \eqno(6.1)
$$
where $ \ N = 1,2, \ldots, \infty. $

 Let $ v = v(n) $ be some non-decreasing positive deterministic function,
$ v(n) \to \infty $ as $ n \to \infty. $
 We purpose that for some $ \psi \in \Psi(a,b) $

$$
\sup_n ||S(n)/\sigma(n)||G(\psi) < \infty. \eqno(6.2)
$$

  We will obtain in this section using (6.1) the rearrangement norm
 estimations for the value

$$
\tau =
|| \sup_n \left[ S(n)/( v(n) \ \sigma(n)) \right]||G(\psi_1), \eqno(6.3)
$$

where at $ p > 1 $

$$
 \psi_1(p) = p \ \psi(p) /(p \  - \ 1).
$$

 In the "entropy" and  "generic chaining" terms in the probabilistic case
$ \mu(X) = 1 $ this estimations are obtained in [13], [16], [40].\par

\vspace{2mm}

{\bf Theorem 2.}  Let $ v = v(n) $ be such that

$$
\sum_{n=1}^{\infty} 1/v \left(2^n \right) < \infty. \eqno(6.4)
$$

Then

$$
||\tau||G(\psi_1) \le C \  \sup_n ||S(n)/\sigma(n)||G(\psi). \eqno(6.5).
$$
{\bf Proof.} We intend to use the inequality (3.3), where instead Pizier
assertion we will use the Doob's inequality.\par

 Choosing the {\it partition } over the  closed
intervals $ W = \{ [A(k), A(k+1) - 1] \} = \{ [A(k), B(k) ] = \{Q(k) \} \} $
of a view:
$$
Q(k) = [A(k), B(k)] = [2^{k-1}, 2^k \ - \ 1], \ k = 1,2,\ldots.
$$
 Suppose for simplicity
$$
\sup_n ||S(n)/\sigma(n)||G(\psi) = 1.
$$
 Let us denote
$$
\tau(k) = \max_{ m \in Q(k)} | S(m)/ (\sigma(m) \ v(m)|;
$$
then

$$
| \tau|_p \le \sum_k |\tau(k)|_p.
$$

Further,

$$
|\tau(k)|_p = \left|\max_{m \in Q(k) } \frac{|S(m)|}{\sigma(m)
 \ v(m)} \right|_p \le
| \max_{m \in Q(k)} |S(m)|/( v(A(k)) \ \sigma(A(k))) ) |_p \le
$$
	
$$
\frac{p}{p \ - \ 1} \cdot \frac{ |S(B(k))|_p }{ v(A(k)) \ \sigma(A(k))} \le
\frac{p}{p \ - \ 1} \cdot \frac{ \psi(p) \ \sigma(B(k))}{ v(A(k)) \ \sigma(A(k))} \le
$$

$$
C_2 \ \psi_1(p) \ 2^{-k \gamma},
$$

where

$$
C_2 = \sup_{n} L(2n)/L(n) < \infty.
$$

 The proposition of theorem 2 follows after summing over $ k. $ \par

For example, if in addition  for $ n \ge 16 $ and for some
$ \Delta = const > 0 $

$$
v(n)  \ge (\log n) \ (\log \log n)^{1 + \Delta},
$$

then

$$
|| \sup_n [ S(n)/( \sigma(n) \ v(n)) ] ||G(\psi_1) \le C \sup_n
|| \ S(n)/\sigma(n) \ || G(\psi) \cdot (1/\Delta). \eqno(6.6)
$$

{\bf Remark 2.} In the "probabilistic" case $ \mu(X) = 1 $ or, equally,
$ \mu(X) < \infty $ the "true" norming function
is $ v(n)  = (\log \log n)^{1/2} $ for martingales with independent
increments, or, in more general case
$ v(n)  = (\log \log n)^{r/2}, \ r = 1,2,\ldots; $  see [13], [16], [40].
This is an open question: what is the true norming function $ v = v(n) $
in the "unbounded" case $ \mu(X) = \infty \ ? $\par

\vspace{3mm}

{\bf 6. Applications into the theory of Fourier series.}\par

\vspace{3mm}

 In this section we intend to obtain the uniform $ G(\psi) $ bounds for
{\it maximal function } for the partial sums of Fourier series. \par

Let $ X = [-\pi,\pi], \ \mu(dx) = dx, \ c(n) = c(n,f) = $
$$
 \int_{-\pi}^{\pi} \exp(i n x) f(x) dx, n = 0,\pm 1,\pm 2 \ldots; \
 2 \pi s_M[f](x)  =
$$
$$
\sum_{ \{n: |n| \le M \} } c(n) \exp(-inx), \ s^*[f] = \sup_{M \ge 1}
|s_M[f]|;
$$
i.e. in the considered case  $ T = \{1,2,3, \ldots. \} $ \par
 Let for some function $ \psi \in \Psi \ f(\cdot) \in G(\psi) $ and
denote  for the values $ p > 1 $

$$
\psi_2(p) = p^4 \ \psi(p)/(p \ - 1 \ )^2.
$$

{\bf Theorem 3.}

$$
||s^*[f]||G(\psi_2) \le C \ ||f||G(\psi). \eqno(7.1)
$$

{\bf Proof} is at the same as in section 6; we use at the same partition

$ W = \{ [A(k), A(k+1) \ - \ 1] \} = \{ [A(k), B(k) ] = \{Q(k) \} \} $
of a view:
$$
Q(k) = [A(k), B(k)] = [2^{k \ - \ 1}, 2^k \ - \ 1], \ k = 1,2,\ldots;
$$

but instead the Doob's inequality we use the following estimation:
at $ p > 1 $

$$
|s^*[f]|_p \le C \  p^4 \ |f|_p/(p \ - 1 \ )^2;
$$
see, for example, [19], p. 183.\par
 The case of Fourier transform (instead Fourier series), the case of wavelet
or Haar's series and multidimensional case $ X = [-\pi, \pi]^d, \ X = R^d,
\ d \ge 2 $ may be considered analogously. See, e.g. [23]. \par


\vspace{5mm}

\begin{center}

{\bf REFERENCES }\\

\end{center}
  1. {\sc Fernique X.} (1975). {\it Regularite des trajectoires des
    function aleatiores gaussiennes.} Ecole de Probablite de
    Saint-Flour, IV – 1974, Lecture Notes in Mathematic. {\bf 480} 1 – 96,    Springer Verlag, Berlin.\\

  2. {\sc Kozachenko Yu. V., Ostrovsky E.I.} (1985). {\it The Banach Spaces
    of random Variables of subgaussian type.}  Theory of Probab. and Math.
      Stat. (in Russian). Kiev, KSU, {\bf 32}, 43 - 57.\\

   3. {\sc Ledoux M., Talagrand M.} (1991) {\it \ Probability in Banach Spaces.}
      Springer, Berlin, MR 1102015.\\

  4. {\sc Ostrovsky E.I.} (1999). {\it Exponential estimations for Random
    Fields and its applications.} (in Russian). Russia, OINPE.\\

  5. {\sc Ostrovsky E.I.} (2002).{\it Exact exponential estimations for random
     field maximum distribution.} Theory Probab. Appl. {\bf 45} v.3,
      281 - 286. \\

   6. {\sc Talagrand M.} (1996). {\it Majorizing measure: The generic chaining.}
       Ann. Probab. {\bf 24} 1049 - 1103. MR1825156 \\

   7. {\sc Talagrand M.} (2001). {\it Majorizing Measures without Measures.}
     Ann. Probab. 29, 411-417. MR1825156 \\

   8. {\sc Talagrand M.} (2005). {\it The Generic Chaining. Upper and
     Lower Bounds of Stochastic Processes.} Springer, Berlin. MR2133757.\\

   9. {\sc Talagrand M.}(1990). {\it Sample boundedness of stochastic
     processes under increment conditions.} Ann. Probab. {\bf 18}, 1 - 49.\\

   10. {\sc Pizier G.} {\it Condition  $ d^/entropic $ assupant la continuite
   de certain processus et applications a $ l^/analyse $ harmonique.} Seminaire $ d^/ analyse $   fonctionnalle. (1980) Exp. 13 p. 23 - 24.  \\

   11. {\sc Karadzhov G.E., Milman M.}{\it Extrapolation theory: new results and
   applications.} J. Approx. Theory, 113 (2005), 38 - 99. \\

   12. {\sc Jawerth B,  Milman M.} {\it Extrapolation theory with applications.}
   Mem. Amer. Math. Soc. 440 (1991).\\

    13. {\sc Ostrovsky E., Sirota L.} {\it Exponential Bounds in the Law of iterated Logarithm for Martingales. } Electronic publications, arXiv:0801.2125v1 [math.PR] 14 Jan 2008.\\

  14. {\sc A.Fiorenza.} {\it Duality and reflexivity in grand Lebesgue spaces.}
  Collectanea Mathematica (electronic version), {\bf 51}, 2, (2000), 131 - 148.\\

  15. {\sc A. Fiorenza and G.E. Karadzhov.} {\it Grand and small Lebesgue spaces and their analogs.} Consiglio Nationale Delle Ricerche, Instituto per le
  Applicazioni del Calcoto Mauro Picine", Sezione di Napoli, Rapporto tecnico n.
  272/03, (2005).\\

 16. {\sc P. Hall, C.C.Heyde.} {\it Martingale Limit Theorems and its Applications.} USA, New York, Academic Press Inc., (1980); \\

  17. {\sc T.Iwaniec and C. Sbordone.} {\it On the integrability of the Jacobian under minimal hypotheses.} Arch. Rat.Mech. Anal., 119, (1992), 129 – 143.\\

  18. {\sc T.Iwaniec, P. Koskela and J. Onninen.} {\it Mapping of finite distortion:  Monotonicity and Continuity,} Invent. Math. 144 (2001), 507 - 531.\\

  19. {\sc Juan Arias de Reyna. } {\it Pointwiese Convergence Fourier Series.}
   New York, Lect. Notes in Math., (2004); \\

  20. {\sc M.A.Krasnoselsky, Ya.B.Rutisky.} {\it Convex functions and Orlicz's
  Spaces.} P. Noordhoff LTD, The Netherland, Groningen, 1961. \\

21. {\sc E.Ostrovsky.} {\it Exponential Orlicz's spaces: new norms and applications.} Electronic Publications, arXiv/FA/0406534, v.1,  (25.06.2004.)\\

22. {\sc E.Ostrovsky, L.Sirota.} {\it Some new rearrangement invariant spaces: theory and applications.} Electronoc publications: arXiv:math.FA/0605732 v1, 29, (May 2006);\\

 23. {\sc E.Ostrovsky, L.Sirota.} {\it Fourier Transforms in Exponential Rearrangement Invariant Spaces.} Electronoc publications: arXiv:math.FA/040639, v1, (20.6.2004.)\\

24. {\sc M.M.Rao. } {\it Measure Theory and Integration.} Basel - New York, John Wiley, Marcel Decker, second Edition, (2004);\\

25. {\sc  M.M. Rao, Z.D.Ren.} {\it Theory of Orlicz Spaces.}  Basel - New York,
Marcel Decker, (1991);\\

26. {\sc M.M. Rao, Z.D.Ren.} {\it Application of Orlicz Spaces.}  Basel - New York, Marcel Decker, (2002);\\

27. {\sc E. Seneta E.} {\it Regularly Varying Functions}. Mir, Moscow edition,
(1985);\\

28. {\sc  Ostrovsky E., Sirota L.} {\it Moment Banach spaces: Theory and
applications.} HAIT Journal of Science and Engineering C. V. 4, Issue 1 - 2,
pp. 233 - 262, (2007). \\

29. {\sc Ostrovsky E., Sirota L.} {\it Nikol'skii-type inequalities in some
rearrangement invariant spaces.} Electronic publications, arXiv 0804.2311v1
 [math.FA], 15 Apr. (2008). \\

  30. {\sc Bennet C., Sharpley R.} {\it Interpolation of operators.} Orlando, Academic Press  Inc., (1988).\\

 31. {\sc Ostrovsky E.I.} (2002). {\it Exact exponential estimations for
   random  field maximum distribution.} Theory Probab. Appl. {\bf 45} v.3,
      281 - 286. \\

32. {\sc Lukomsky S.F.} {\it About convergence of Walsh series in the spaces    nearest to } $ L_{\infty}. $ Matem. Zametky, 2001, v.20 B.6,p. 882 -   889.(Russian).\\

   33. {\sc Astashkin S.V.} {\it About interpolation spaces of sum spaces, generated by Rademacher system. } RAEN, issue MMMIU, 1997, v.1 $ N^o $ 1, p.
8-35.\\

   34. {\sc  Capone C., Fiorenza A., Krbec M.} {\it On the Extrapolation Blowups in the  $ L_p $ Scale. }\\

   35. {\sc Astashkin S.V.} {\it Some new Extrapolation Estimates for the Scale
   of   $ L_p \ - $  Spaces.} Funct. Anal. and Its Appl., v. 37 $ N^o $ 3 (2003),  73 - 77. \\

    36. {\sc Dudley R.M.} {\it The sizes of compact of Hilbert space and continuity  of Gaussian processes.} J. Functional Analysis. (1967) B. 1 pp. 290 - 330. \\

    37. {\sc Musielac J.} {\it Orlicz Spaces  and Modular Spaces.} Springer Verlag. 2002.\\

  38. {\sc Harjulehto P., Hanstz P. and Pere M.}{\it Variable exponent
Sobolev Spaces.} Funct. Approx. Comment. Math., {\bf 36} (2006), 79 \ - \ 94.\\

 39. {\sc Harjulehto P., Hanstz P. and Pere M.} {\it Variable exponent
Lebesque spaces and Hardy-Littlewood maximal operator.} Real Anal. Exchange, {\bf 30}(2004), Preprint.  \\

  40. {\sc Ostrovsky E.I.} {\it Exponential Bounds  in the Law of Iterated
Logarithm in Banach Space.} (1994), Math. Notes, {\bf 56}, 5, p. 98 - 107. \\

 41. {\sc Ostrovsky E., Rogover E.} {\it Exact exponential bounds for the random
Fields Maximum Distribution via the majoring Measures (generic Chaining).
Electronic Publications, arXiv:o802v1 [math.PR], 4 Feb 2008.\\

42.{\sc Ostrovsky E., Sirota L.} {\it Moment Banach Spaces: Theory and Applications.} HAIT Journal of Science and Engineering, C,  V. 4 Issues
1  - 2, pp. 233 – 262.\\


\begin{center}

 \hspace{50mm} {\sc Department of Mathematics}\\

  \hspace{50mm} {\sc Bar-Ilan University}\\

  \hspace{50mm} {\sc Ramat Gan, Ben Gurion street, 2}\\

  \hspace{50mm} {\sc Israel \ 76521}\\

  \hspace{50mm} {\sc E-Mail}: \ galo@list.ru \\

  \hspace{50mm} {\sc E-Mail}: \ rogovee@gmail.com \\

\vspace{4mm}

{\bf Ostrovsky E.}\\
\vspace{4mm}

 Address: Ostrovsky E., ISRAEL, 76521, Rehovot, \ Shkolnik street.
5/8. Tel. (972)-8- 945-16-13.\\
\vspace{3mm}
e - mail: {\bf Galo@list.ru}\\

\vspace{6mm}

{\bf Rogover E.}\\
\vspace{3mm}
 Address: Rogover E., ISRAEL, 84105, Ramat Gan.

\vspace{4mm}
e - mail: {\bf rogovee@gmail.com }\\

\vspace{3mm}

{\it Department of Mathematics and Statistics, Bar-Ilan University,
59200, Ramat Gan, Israel.}\\
e-mail: \ galo@list.ru \\

{\it Department of Mathematics and Statistics, Bar-Ilan University,
59200, Ramat Gan, Israel.}\\
e - mail: \ sirota@zahav.net.il \\

\end{center}

\end{document}